\documentclass[DIV=14,letterpaper]{scrartcl}

\usepackage{tikz}
\usetikzlibrary{calc}
\usetikzlibrary{matrix}
\usepackage{amsmath,amssymb,amsfonts,amsthm,subfig}
\usepackage[pdftex]{hyperref} 

\usepackage{graphicx}
\allowdisplaybreaks

\newcommand{\V}[1]{\mbox{\boldmath $ #1 $}}

\newcommand{\bey}{\begin{eqnarray}}
\newcommand{\eey}{\end{eqnarray}}

\newcommand{\beq}{\begin{equation}}
\newcommand{\eeq}{\end{equation}}
\theoremstyle{plain}% default

\theoremstyle{definition}

\theoremstyle{remark}

\title{An Introduction to MMPDElab}

\author{Weizhang~Huang%
\thanks{Department of Mathematics, the University of Kansas, Lawrence, KS 66045, U.S.A. 
({\em whuang@ku.edu}).}
}

\date{}

\begin{document}
\vskip 1cm
\maketitle

% section 1
\section{Introduction}

\textbf{MMPDElab} is a package written in MATLAB\footnote{MATLAB$^{\tiny{\textregistered}}$
is a trademark of The MathWorks, Inc., Natick, MA 01760.}
 for adaptive mesh movement and adaptive moving mesh P1 finite element solution of second-order partial
different equations (PDEs) having continuous solutions in one, two, and three spatial dimensions.
It uses simplicial meshes, i.e., line segments in one dimension, triangles in two dimensions, and tetrahedra
in three dimensions. The adaptive mesh movement is based
on the new implementation \cite{HK2014,HK2015} of the moving mesh partial differential equation (MMPDE)
method \cite{BHR09,CHR99b,HRR94b,HRR94a,HR97b,HR99,HR11}. The mesh equation is integrated using
either \verb|ode45| (an explicit MATLAB ODE solver) or \verb|ode15s| (an implicit MATLAB ODE solver)
while physical PDEs are discretized in space using P1 conforming finite elements on moving meshes
and integrated in time with the fifth-order Radau IIA method (an implicit Runge-Kutta method)
with a two-step error estimator \cite{Montijano2004} for time step selection.
More information on the moving mesh P1 finite element method can be found from
recent applications such as those found in
\cite{DHHLY2018,HKS2015,LHQ2018,NH2017,NH2019,YH2018,ZhangFei2017}. 

The source code of \textbf{MMPDElab} can be downloaded at 
\begin{itemize}
\item https://whuang.ku.edu/MMPDElab/mmpdelabv1.html
\item https://github.com/weizhanghuang/MMPDElab
\end{itemize}

\vspace{20pt}

The functions in MMPDElab can be grouped into three categories:
\begin{itemize}
\item Matrix operations (with names in the form \verb|Matrix_xxx|)
\item Mesh movement (with names in the form \verb|MovMesh_xxx|)
\item Moving mesh P1 finite element solution (with names in the form \verb|MovFEM_xxx|)
\end{itemize}

The functions in the first category \verb|Matrix_xxx| perform vectorized computation of basic
matrix operations such as multiplication, inversion, and finding transposes and determinants
for arrays of matrices of small size (typically $3\times 3$ or smaller). These operations are
used by functions in the other two categories which will be explained in the subsequent sections.

We now introduce notation whose understanding is crucial to the use of the package.
A simplicial mesh or a simplicial triangulation, $\mathcal{T}_h$,
of $N$ elements and $N_v$ vertices in $d$-dimensions ($d = 1$, 2, or 3) is represented in MATLAB by the matrices
$X$ and $tri$, where $X$ is a matrix of size $N_v \times d$ containing the coordinates of the vertices and
$tri$ is a matrix of size $N \times (d+1)$ listing the connectivity of the mesh. More specifically, $X(i,:)$ gives
the coordinates of the $i$th vertex $\V{x}_i$ while $tri(j,:)$ contains the global IDs of the vertices of the $j$th element.
In \textbf{MMPDElab}, $npde$ components of the physical solution at the vertices
are given by the matrix $u$ of size $N_v \times npde$, i.e., $u(i, :)$ contains the values of $u$ at the $i$th vertex.
Its derivatives with respect to the physical coordinate $\V{x} = (x_1, x_2, ..., x_d)^T$ are saved in the form
\[
du = \left [ (\nabla u^{(1)})^T, ..., (\nabla u^{(npde)})^T\right ]_{N_v \times (d \ast npde)},
\]
where $u^{(k)}$ ($k = 1, ..., npde$) is the $k$th component of $u$ and $\nabla$ is the gradient operator.
Thus,
\[
du(i,:) = \left [ \frac{\partial u^{(1)}}{\partial x_1}, ..., \frac{\partial u^{(1)}}{\partial x_d}, ..., 
\frac{\partial u^{(npde)}}{\partial x_1}, ..., \frac{\partial u^{(npde)}}{\partial x_d} \right ](\V{x}_i),
\quad i = 1, ..., N_v.
\]
The metric tensor or the monitor function, $\mathbb{M}$, is calculated at the vertices and
saved in the form
\[
M(i,:) = \left [\mathbb{M}_{11}, ..., \mathbb{M}_{d1}, ...,  \mathbb{M}_{1d}, ..., \mathbb{M}_{dd}\right ](\V{x}_i),
\quad i = 1, ..., N_v.
\]
That is, $M$ has the size $N_v \times (d \ast d)$, with each row containing the entries of a matrix of size $d\times d$.
$M$ is a good example of an array of matrices of small size. It is emphasized that when a moving mesh function
is called, the mesh connectivity is kept fixed while the location of the vertices varies. The user can
decide whether or not to change the connectivity in between the calls.

To conclude this section, I am deeply thankful for many colleagues and former graduate students
for their invaluable discussion and comments. I am particularly grateful to Dr. Lennard Kamenski
who was involved in the project at the early stage. 

\vspace{10pt}

\begin{verbatim}
    MMPDElab is a package written in MATLAB  for adaptive mesh 
    movement and adaptive moving mesh P1 finite element solution 
    of partial different equations having continuous solutions.
    
    Copyright (C) 2019  Weizhang Huang (whuang@ku.edu)

    MMPDElab is free software: you can redistribute it and/or modify
    it under the terms of the GNU Affero General Public License as
    published by the Free Software Foundation, either version 3 of the
    License, or (at your option) any later version.

    MMPDElab is distributed in the hope that it will be useful,
    but WITHOUT ANY WARRANTY; without even the implied warranty of
    MERCHANTABILITY or FITNESS FOR A PARTICULAR PURPOSE.
    See the GNU Affero General Public License at
    <https://www.gnu.org/licenses/>.
\end{verbatim}

% section 2
\section{Adaptive mesh movement}

The adaptive mesh movement can be carried out by calling
\verb|MovMesh()| (based on the $\xi$-formulation of the MMPDE moving mesh method \cite{HK2014,HK2015}),
\verb|MovMesh_XM()| (based on the $x$-formulation of the MMPDE moving mesh method),
or \verb|MovMesh_X()| (based on the $x$-formulation of the MMPDE moving mesh method
with the metric tensor $\mathbb{M} = I$, i.e., without mesh adaptation). The corresponding MMPDE
is defined as a gradient flow equation of the meshing functional developed in \cite{Hua01b}
based on mesh equidistribution and alignment (with its parameters being chosen as
$p = 1.5$ and $\theta = 1/3$). The headers of these functions read as

\begin{verbatim}
[Xnew,Ih,Kmin] = MovMesh(tspan,Xi_ref,X,M,tau,tri,tri_bf,nodes_fixed, ...
                         mmpde_int_method,dt0,abstol)
                                                          
[Xnew,Ih,Kmin] = MovMesh_XM(tspan,X,M,tau,tri,tri_bf,nodes_fixed, ...
                            mmpde_int_method,dt0,abstol,Xi_ref)

[Xnew,Ih,Kmin] = MovMesh_X(tspan,X,tau,tri,tri_bf,nodes_fixed, ...
                           mmpde_int_method,dt0,abstol,Xi_ref)
\end{verbatim}

These functions integrate the corresponding moving mesh equations over a time period
specified by {\em tspan}. All of the meshes, {\em X} (the current mesh),
{\em Xnew} (the new mesh), and {\em Xi\_ref} (the reference computational mesh),
are assumed to have the same number of vertices and elements and the same
connectivity (specified by {\em tri}). The input and output variables are explained in the following.
\begin{itemize}
\item {\em tspan} is a vector specifying the time interval for integration.
\item {\em X}, of size $N_v \times d$, contains the coordinates of vertices of the current mesh.
\item {\em Xi\_ref}, of size $N_v \times d$, contains the coordinates of vertices of the reference computational mesh.
	This mesh, typically chosen as the initial physical mesh,
	is a mandatory input for \verb|MovMesh()| but is optional for \verb|MovMesh_XM()|
	and \verb|MovMesh_X()|. In the latter case, when {\em Xi\_ref} is not supplied,
	the uniformity of the new physical mesh measured in the metric $\mathbb{M}$ is made with
	reference to an equilateral simplex.
\item {\em M}, of size $N_v \times (d\ast d)$, contains the values of the metric tensor $\mathbb{M}$
        at the vertices of {\em X}. More specifically, {\em M(i,1:d$\ast$d)} gives the metric
	tensor at the $i$th vertex, i.e.,
	$ [\mathbb{M}_{11}, ..., \mathbb{M}_{d1}, ..., \mathbb{M}_{1d}, ..., \mathbb{M}_{dd}](\V{x}_i)$.
\item {\em tau} is the positive parameter used for adjusting the time scale of mesh movement.
\item {\em tri}, of size $N \times (d+1)$, lists the connectivity for all meshes.
\item {\em tri\_bf}, of size $N_{bf} \times d$, specifies the boundary facets for all meshes, with each row containing
	the IDs of the vertices of a facet on the boundary. A boundary
	facet consists of a point in 1D, a line segment (with two vertices) in 2D, or a triangle (with three
	vertices) in 3D. {\em tri\_bf} can be computed using the Matlab function {\em freeBoundary} in 2D and 3D.
\item {\em nodes\_fixed} is a vector containing the IDs of the vertices, such as corners, which are not allowed to move
	during the mesh movement.
\item {\em mmpde\_int\_method} is an optional input variable, specifying that either \verb|ode15s| (implicit) or
	\verb|ode45| (explicit) is used to integrate the moving mesh equation. The default is \verb|ode15s|.
\item {\em dt0} is an optional input variable specifying the initial time stepsize that is used in
	the time integration of the mesh equation. The default is  {\em dt0} = ({\em tspan}(end)-{\em tspan}(1))/10.
\item {\em abstol} is an optional input variable specifying the absolute tolerance used for time step selection
	in the time integration of the mesh equation. The default is {\em asbstol} = 1e-6 for \verb|ode15s|
	and 1e-8 for \verb|ode45|.
\item {\em Xnew}, of size $N_v \times d$, contains the coordinates of vertices of the new mesh.
\item {\em Ih} is an optional output variable giving the value of the meshing functional at the new mesh.
\item {\em Kmin} is an optional output variable giving the minimal element volume.
\end{itemize}

In addition to \verb|MovMesh()|, \verb|MovMesh_XM()|, and \verb|MovMesh_X()|,
the following functions can also be used by the user.
\begin{enumerate}
\item \verb|[X,tri] = MovMesh_circle2tri(jmax) |
	This function creates a triangular mesh ({\em X}, {\em tri}) for the unit circle.
\item \verb|[X,tri] = MovMesh_cube2tet(x,y,z) |
	This function creates a tetrahedral mesh ({\em X}, {\em tri}) from the cuboid mesh specified by {\em x},
	{\em y}, and {\em z} for a cuboid domain. Each subcuboid is divided into 6 tetrahedra.
\item \verb|V = MovMesh_freeBoundary_faceNormal(X,tri,tri_bf) |
	This function computes the unit outward normals for the boundary facets.
	{\em V} has the size of $N_{bf} \times d$.
\item \verb|V = MovMesh_freeBoundary_vertexNormal(X,tri,tri_bf) |
	This function computes the unit outward normals for the boundary vertices.
	{\em V} has the size of $N_{v} \times d$, with the normals for the interior vertices being set to be
	$[1, ..., 1]^T/\sqrt{d}$.
\item \verb|[Grad,Hessian] = MovMesh_GradHessianRecovery(u,X,tri,tri_bf) |
	This function computes the gradient and Hessian of function $u$ at the vertices using
	centroid-vortex-centroid-vertex volume-weighted average.
\item \verb|Grad = MovMesh_GradKRecovery(u,X,tri,tri_bf) |
	This function computes the gradient of function $u$ on the elements.
\item \verb|Grad = MovMesh_GradRecovery(u,X,tri,tri_bf) |
	This function computes the gradient of function $u$ at the vertices using volume averaging.
\item \verb|fnew = MovMesh_LinInterp(f,X,QP,tri,tri_bf,useDelaunayTri) |
	This function performs linear interpolation of {\em f} (defined on {\em X}) at query points QP using
	triangulation or Delaunay triangulation. {\em useDelaunayTri} is a logical variable with
	value {\em true} or {\em false}.
\item \verb|[X,tri] = MovMesh_MeshMerge(X1,tri1,X2,tri2) |
	This function merges two non-overlapping meshes ({\em X1,tri1}) and ({\em X2,tri2}) which
	may or may not have common boundary segments.
\item \verb|[Qgeo,Qeq,Qali] = MovMesh_MeshQualMeasure(X,tri,M,Linf_norm,Xi_ref) |
	This function computes the geometric, equidistribution, and alignment measures
	(in maximum norm or $L^2$ norm in $\xi$) for the mesh ({\em X, tri})
	according to the metric tensor. Here, both {\em Linf\_norm} and {\em Xi\_ref} are optional input variables.
\item \verb|[Qmax,Ql2] = MovMesh_MeshQualMeasure2(X,tri,M,Xi_ref) |
	This function computes the maximum and $L^2$ norm of the mesh quality measure based on
	a single condition combining both equidistribution and alignment.
	{\em Xi\_ref} is an optional input variable.
\item \verb|[X,tri] = MovMesh_MeshRemoveNodes(X1,tri1,ID) |
	This function removes the nodes listed in {\em ID} from the existing mesh ({\em X1,tri1}).
\item \verb|[XF,TriF,TriF_parent] = MovMesh_MeshUniformRefine(X,Tri,Level) |
	This function uniformly  refines a simplicial mesh ({\em Level}) times or levels. On each level,
	an element is refined into $2^d$ elements.
\item \verb|M = MovMesh_metric_arclength(u,X,tri,tri_bf) |
	This function computes the arclength metric tensor.
\item \verb|MC = MovMesh_metric_F2C(M,Tri,Tri_parent,TriC) |
	This function computes the metric tensor on a coarse mesh from the metric tensor
	defined on a fine mesh.
\item \verb|M = MovMesh_metric_intersection(M1,M2) |
	This function computes the intersection of two symmetric and positive definite matrices.
	When \verb|M1| and \verb|M2| are diagonal, i.e., $\verb|M1|  = \text{diag}(\alpha_1, ..., \alpha_d)$
	and $\verb|M2|  = \text{diag}(\beta_1, ..., \beta_d)$, then
	$\verb|M|  = \text{diag}(\max(\alpha_1, \beta_1), ..., \max(\alpha_d, \beta_d))$.
	The intersection of two general symmetric and positive definite matrices is defined
	similarly through simultaneous diagonalization.
\item \verb|M = MovMesh_metric_iso(u,X,tri,tri_bf,alpha,m) |
	This function computes the isotropic metric tensor based on the $L^2$ norm or
	the $H^1$ seminorm of linear interpolation error ($m = 0$ or $m = 1$).
\item \verb|MM = MovMesh_metric_smoothing(M,ncycles,X,tri) |
	This function smooths the metric tensor {\em ncycles} times by local averaging.
\item \verb|M = MovMesh_metric(u,X,tri,tri_bf,alpha,m) |
	This function computes the metric tensor based on the $L^2$ norm or
	the $H^1$ seminorm of linear interpolation error ($m = 0$ or $m = 1$).
\item \verb|[X,tri] = MovMesh_rect2tri(x,y,job) |
	This function creates a triangular mesh ({\em X}, {\em tri}) from the rectangular mesh
	specified by {\em x} and {\em y} for a rectangular domain. Each rectangle is divided
	into 2 (for {\em job} = 2 or 3) or 4 (for {\em job} = 1) triangles.
\item \verb|M1 = Matrix_ceil(M,beta) |
	This function puts a ceiling on the eigenvalues of symmetric and positive definite matrix
	{\em M} such that $\lambda_{max}(M1) \le \beta$.
\end{enumerate}

Examples using these functions include \verb|ex1d_1.m|, \verb|ex2d_1.m|, \verb|ex2d_2_Lshape.m|,
\verb|ex2d_3_hole.m|, \verb|ex2d_4_horseshoe.m|, and \verb|ex3d_1.m| in the subdirectory \verb|./examples|.

\vspace{20pt}

\textbf{Troubleshooting.}
Occasionally one may see an error message like
{\small
\begin{verbatim}
Error using triangulation
The coordinates of the input points must be finite values; Inf and NaN are not permitted.

Error in MovMesh>MovMesh_rhs (line 296)
         TR = triangulation(tri2,XI2);
\end{verbatim}
}
\noindent when calling \verb|MovMesh()|, \verb|MovMesh_XM()|, or \verb|MovMesh_X()|. Typically this is caused
by a stability issue when integrating the MMPDE, and using a smaller initial time step {\em dt0} (e.g.,
{\em dt0 = 1e-6}) may solve the problem.

\vspace{20pt}

\textbf{The generation of initial meshes.} This package includes a few functions for generating
initial meshes  for simple domains, such as
\verb| MovMesh_circle2tri()|, \verb|MovMesh_cube2tet()|, and\\ \verb|MovMesh_rect2tri()|.
For complex domains, one can use MATLAB function \verb|delaunayTriangulation()| or
other automatic mesh generators such as DistMesh (in MATLAB) \cite{Persson2012}
and TetGen (in C++) \cite{HangSi}.

% section 3
\section{Adaptive mesh movement P1 finite element solution of PDEs}

This package aims to solve the system of PDEs in the weak form: find $u = [u^{(1)}, ..., u^{(npde)}]
\in H^1(\Omega) \otimes \cdots \otimes H^1(\Omega)$ such that
\begin{align}
\label{PDE-1}
&\sum_{i=1}^{npde} \int_\Omega F_i(\nabla u, u, u_t, \nabla v^{(i)}, v^{(i)}, \V{x}, t) d \V{x}
+ \sum_{i=1}^{npde} \int_{\Gamma_N^{(i)}} G_i(\nabla u, u, v^{(i)}, \V{x}, t) d s = 0, \\
& \qquad \qquad \qquad \qquad \qquad \qquad \qquad \qquad \qquad \quad \quad
\forall v^{(i)} \in V^{(i)}, \quad i = 1, ..., npde, \quad 0 < t \le T
\notag
\end{align}
subject to the Dirichlet boundary conditions
\begin{align}
R_i(u, \V{x}, t) = 0, \qquad \text{ on } \Gamma_D^{(i)}, \quad  i=1, ..., npde
\label{BC-1}
\end{align}
where for $i= 1, ..., npde$, $\Gamma_D^{(i)}$ and $\Gamma_N^{(i)}$ are the boundary segments
corresponding to the Dirichlet and Neumann boundary conditions for $u^{(i)}$,  respectively,
$\Gamma_D^{(i)}  \cup \Gamma_N^{(i)} = \partial \Omega$, and
$V^{(i)} = \{ w \in H^1(\Omega), \; w = 0 \text{ on } \Gamma_D^{(i)} \}$.
The headers of \verb|MovFEM()| (Initial-Boundary-Value-Problem solver) and
\verb| MovFEM_bvp()| (Boundary-Value-Problem solver)  read as
\begin{verbatim}
[Unew,dt0,dt1] = MovFEM(t,dt,U,X,Xdot,tri,tri_bf,pdedef, ...
                        fixed_step,relTol,absTol,direct_ls,ControlWeights)
                      
Unew = MovFEM_bvp(U,X,tri,tri_bf,pdedef,nonlinearsolver,MaxIter,Tol)
\end{verbatim}

\verb|MovFEM()| integrates the system of PDEs on a moving mesh over a time step. Its input and output variables
are explained in the following.
\begin{itemize}
\item {\em t} is the  current time.
\item {\em dt} is the intended time stepsize for integrating the physical PDEs.
\item {\em U}, of size $N_v \times npde$, is the current solution.
\item {\em X}, of size $N_v \times d$,  contains the coordinates of vertices of the current mesh.
\item {\em Xdot}, of size $N_v \times d$, is the nodal mesh velocity.
\item {\em tri}, of size $N \times (d+1)$, lists the connectivity for all meshes.
\item {\em tri\_bf}, of size $N_{bf} \times d$, specifies the boundary facets for all meshes.
\item {\em pdedef} is a structure used to define the PDE system in the weak form. It has 5 fields.
	\begin{enumerate}
	\item[(i)] {\em pdedef.bfMark}, of size $N_{bf} \times 1$, is used to mark the boundary segments
		(boundary facets). This marking is passed to the definitions of boundary integrals and
		Dirichlet boundary conditions.
	\item[(ii)]	{\em pdedef.bftype},  of size $N_{bf} \times npde$, specifies the types of boundary condition
		on boundary facets whose numbering is based on {\em tri\_bf}. {\em pdedef.bftype} = 0 for Neumann
		BCs and {\em pdedef.bftype} = 1 for Dirichlet BCs. For example, {\em pdedef.bftype(3,2)} = 1
		means that variable $u^{(2)}$ has a Dirichlet BC on the 3rd boundary facet
		while {\em pdedef.bftype(2,1) = 0} specifies that variable $u^{(1)}$ has a Neumann BC
		on the 2nd boundary facet.
	\item[(iii)] \verb|F = pdedef.volumeInt(du, u, ut, dv, v, x, t, i) | This function is used to define
		$F_i$ in the weak form (\ref{PDE-1}), where $v$ and $dv$ are the test function $v^{(i)}$ and its gradient.
	\item[(iv)] \verb|G = pdedef.boundaryInt(du, u, v, x, t, i, bfMark) | This function is used to define
		 $G_i$ in the weak form (\ref{PDE-1}), where $v$ is the test function $v^{(i)}$.
	\item[(v)] \verb|R = pdedef.dirichletRes(u, x, t, i, bfMark) | This function is used to define $R_i$
		in (\ref{BC-1}).
	\end{enumerate}
\item {\em fixed\_step} is an optional input logical variable, indicating whether or not a fixed
	step is used in time integration. The default is {\em false}.
\item {\em relTol} and {\em absTol} are optional input variables for the relative and absolute
	tolerances used for time step selection.
	The defaults are {\em relTol} = 1e-4 and {\em absTol} = 1e-6.
\item {\em direct\_ls} is an optional input logical variable, indicating whether or not the direct sparse
	matrix solver is used for solving linear algebraic systems. 
	When {\em direct\_ls} = {\em false}, the BiConjugate Gradients Stabilized Method
	\verb|bicgstab| is used. The default is {\em true}.
\item {\em ControlWeights} is an optional input variable which is nonnegative vector of size
	$(N_v\ast npde) \times 1$ used to define the weights of the components of the solution
	for the error estimation used in time step selection.
\item {\em Unew}, of size $N_v \times npde$,  is the new solution at time {\em t + dt0}.
\item {\em dt0} is the time stepsize actually used to integrate the physical PDEs.
\item {\em dt1} is the time step size predicted for the next step.
\end{itemize}

The input and output variables for \verb|MovFEM_bvp()| are similar to those of \verb|MovFEM()|.
The same weak form (\ref{PDE-1}) and (\ref{BC-1}) is used for both IBVPs and BVPs.
In the latter case, $t$ is a parameter that is not used. Here we list the variables used
only in the BVP solver.
\begin{itemize}
\item {\em nonlinearsolver} is an optional input variable for the method used for solving nonlinear
	algebraic systems, with the choices being {\em newtons} and  \verb|fsolve|.
	The defacult is  \verb|fsolve|.
\item {\em MaxIter} is an optional input variable for the maximum number of iterations allowed
	for the solution of nonlinear algebraic systems. The default is {\em MaxIter} = 300.
\item {\em Tol} is an  optional input variable for the tolerance used in the solution of nonlinear algebraic
	systems. The default is {\em Tol} = 1e-6.
\end{itemize}

The following two functions are used to compute the error when the exact solution
is available in the form \verb|U = uexact(t,x,varargin)|.
\begin{enumerate}
\item \verb|err = MovFEM_Error_P1L2(uexact,t,X,U,tri,tri_bf,varargin) |
	This function computes the $L^2$ norm of the error in P1 finite element approximation.
\item \verb|err = MovFEM_Error_P1Linf(uexact,t,X,U,tri,tri_bf,varargin) |
	This function computes the $L^\infty$ norm of the error in P1 finite element approximation.
\end{enumerate}

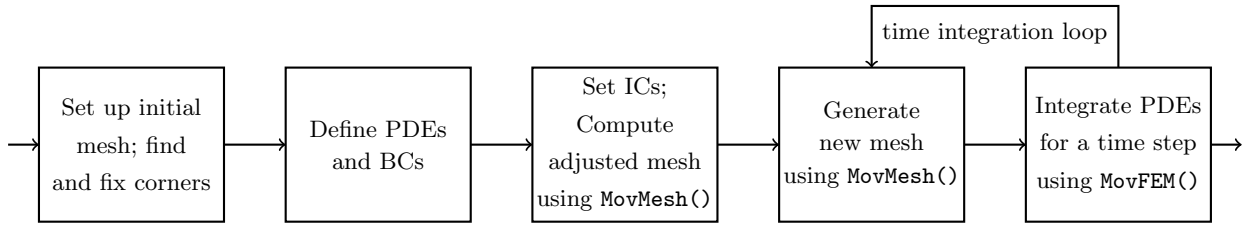
\begin{figure}[h]
\centering
\tikzset{my node/.code=\ifpgfmatrix\else\tikzset{matrix of nodes}\fi}
{\footnotesize
\begin{tikzpicture}[every node/.style={my node},scale=0.41]
\draw[->,thick] (-1,2.5)--(0,2.5);
\draw[thick] (0,0) rectangle (6,5);
\node (node1) at (3,2.5) {Set up initial\\ mesh; find\\ and fix corners\\};
\draw[->,thick] (6,2.5)--(8,2.5);
\draw[thick] (8,0) rectangle (14,5);
\node (node2) at (11,2.5) {Define PDEs\\ and BCs\\};
\draw[->,thick] (14,2.5)--(16,2.5);
\draw[thick] (16,0) rectangle (22,5);
\node (node3) at (19,2.5) {Set ICs; \\Compute\\ adjusted mesh\\ using \verb|MovMesh()|\\};
\draw[->,thick] (22,2.5)--(24,2.5);
\draw[thick] (24,0) rectangle (30,5);
\node (node4) at (27,2.5) {Generate\\ new mesh\\ using \verb|MovMesh()|\\};
\draw[->,thick] (30,2.5)--(32,2.5);
\draw[thick] (32,0) rectangle (38,5);
\node (node4) at (35,2.5) {Integrate PDEs\\ for a time step\\ using \verb|MovFEM()|\\};
\draw[->,thick] (38,2.5)--(39,2.5);

\draw[->,thick] (35,5)--(35,7)--(27,7)--(27,5);
\node[below] at (31,7) {time integration loop\\};

\end{tikzpicture}
}
\caption{An MP (Mesh PDE -- Physical PDE) procedure for moving mesh solution of IBVPs.}
\label{fig:ibvp-solver-1}
\end{figure}

In the following we give several examples to explain how to define the PDE system
through {\em pdedef}. More examples can be found in the subdirectory \verb|./examples|.
A typical flow chart for those examples is shown in Fig.~\ref{fig:ibvp-solver-1}.

\vspace{20pt}

% section 3.1
\subsection{Burgers' equation in 1D}
% exam: burgers equation in 1D
%\begin{exam}
\label{exam-bergers1d}
This example, implemented in \verb|ex1d_burgers.m|, is the IBVP of Burgers' equation in 1D,
\begin{equation}
\label{burgers1d-1}
u_t = \epsilon u_{xx} - u u_x, \quad x \in \Omega \equiv (0,1), \quad t \in (0, 1]
\end{equation}
subject to the Dirichlet boundary condition
\begin{equation}
\label{burgers1d-1-BC}
u(t,x) = u_{exact}(t,x), \quad x \text{ on } \partial \Omega, \quad t \in (0, 1]
\end{equation}
and the initial condition
\begin{equation}
\label{burgers1d-1-IC}
u(0,x) =  u_{exact}(0,x), \quad x \in \Omega
\end{equation}
where $\epsilon = 10^{-3}$ and
\begin{equation}
\label{burgers1d-1-exac}
u_{exact}(t, x) = \frac{0.1 e^{\frac{-x+0.5-4.95t}{20 \epsilon}}
+ 0.5 e^{\frac{-x+0.5-0.75t}{4 \epsilon}} + e^{\frac{-x+0.375}{2 \epsilon}} }
{e^{\frac{-x+0.5-4.95t}{20 \epsilon}} + e^{\frac{-x+0.5-0.75t}{4 \epsilon}}
+ e^{\frac{-x+0.375}{2 \epsilon}} } .
\end{equation}
The weak formulation of this example reads as
\begin{equation}
\label{burgers1d-2}
\int_\Omega (u_t v + \epsilon u_x v_x + u u_x v) d x = 0,\quad \forall v \in V \equiv H^1_0(\Omega) .
\end{equation}
The definition of this example in the code is given as
\begin{verbatim}
   % define PDE system and BCs
   
   % all bcs are dirichlet so no need for marking boundary segments
   pdedef.bfMark = ones(Nbf,1);
   pdedef.bftype = ones(Nbf,npde);
   
   pdedef.volumeInt = @pdedef_volumeInt;
   pdedef.boundaryInt = @pdedef_boundaryInt;
   pdedef.dirichletRes = @pdedef_dirichletRes;

... ...

function F = pdedef_volumeInt(du, u, ut, dv, v, x, t, ipde)
   global epsilon;
   F = ut(:,1).*v(:) + epsilon*du(:,1).*dv(:,1) + u(:,1).*du(:,1).*v(:); 
    
function G = pdedef_boundaryInt(du, u, v, x, t, ipde, bfMark)
   G = zeros(size(x,1),1);

function Res = pdedef_dirichletRes(u, x, t, ipde, bfMark)
   Res = u - uexact(t, x);
\end{verbatim}
%\end{exam}

% section 3.2
\subsection{The heat equation in 2D}
% exam: heat equation in 2D
%\begin{exam}
\label{exam-heat2d}
This example, implemented in \verb|ex2d_heat.m|, is the IBVP for the heat equation in 2D,
\begin{equation}
\label{heat2d-1}
u_t = u_{xx} + u_{yy} + (13\pi^2-1) u_{exact}(t,x, y),
\quad (x,y) \in \Omega \equiv (0,1)\times (0, 1), \quad t \in (0, 1]
\end{equation}
subject to the boundary conditions
\begin{equation}
\label{heat2d-1-BC}
\begin{cases}
u(t,x,y) = 0, & \quad (x,y) \text{ on } x = 0 \text{ and } y = 0, \quad t \in (0, 1] \\
\frac{\partial u}{\partial x} = 2\pi e^{-t} \sin(3 \pi y), &\quad (x,y) \text{ on } x = 1, \quad t \in (0, 1] \\
\frac{\partial u}{\partial y} = - 3\pi e^{-t} \sin(2\pi x), &\quad (x,y) \text{ on } y = 1, \quad t \in (0, 1]
\end{cases}
\end{equation}
and the initial condition
\begin{equation}
\label{heat2d-1-IC}
u(0,x,y) =  u_{exact}(0,x, y), \quad (x,y) \in \Omega .
\end{equation}
This problem has the exact solution
\begin{equation}
\label{heat2d-1-exac}
u_{exact}(t,x, y) = e^{-t} \sin(2\pi x) \sin(3 \pi y).
\end{equation}
The weak formulation reads as
\begin{align}
\label{heat2d-2}
& \int_\Omega \left [ (u_t v + u_x v_x + u_y v_y ) - (13\pi^2-1)\, v\, u_{exact}(t,x, y) \right ] d x dy
\\
& \qquad + \int_0^1 \left ( - 2\pi e^{-t} \sin(3 \pi y) \right ) v(1, y) d y 
+ \int_0^1 \left ( 3\pi e^{-t} \sin(2 \pi x) \right ) v(x, 1)  d x = 0,\quad \forall v \in V 
\notag
\end{align}
where $V = \{ w \in H^1(\Omega), w = 0 \text{ on } x = 0 \text{ and } y = 0\}$.
The definition of  this example in the code is given as
\begin{verbatim}
   % define PDE system and BCs
   
   % mark boundary segments
   pdedef.bfMark = ones(Nbf,1); % for y = 0 (b1)
   Xbfm = (X(tri_bf(:,1),:)+X(tri_bf(:,2),:))*0.5;
   pdedef.bfMark(Xbfm(:,1)<1e-8) = 4; % for x = 0 (b4)
   pdedef.bfMark(Xbfm(:,1)>1-1e-8) = 2; % for x = 1 (b2)
   pdedef.bfMark(Xbfm(:,2)>1-1e-8) = 3; % for y = 1 (b3)
   
   % define boundary types
   pdedef.bftype = ones(Nbf,npde);
   % for neumann bcs:
   pdedef.bftype(pdedef.bfMark==2|pdedef.bfMark==3,npde) = 0;
            
   pdedef.volumeInt = @pdedef_volumeInt;
   pdedef.boundaryInt = @pdedef_boundaryInt;
   pdedef.dirichletRes = @pdedef_dirichletRes;   
... ...

function F = pdedef_volumeInt(du, u, ut, dv, v, x, t, ipde)
   F = (13*pi*pi-1)*uexact(t,x);
   F = ut(:,1).*v(:)+du(:,1).*dv(:,1)+du(:,2).*dv(:,2)-F.*v(:); 

function G = pdedef_boundaryInt(du, u, v, x, t, ipde, bfMark)
   G = zeros(size(x,1),1);
   ID = find(bfMark==2);
   G(ID) = -2*pi*exp(-t)*sin(3*pi*x(ID,2)).*v(ID);
   ID = find(bfMark==3);
   G(ID) = 3*pi*exp(-t)*sin(2*pi*x(ID,1)).*v(ID);

function Res = pdedef_dirichletRes(u, x, t, ipde, bfMark)
   Res = zeros(size(x,1),1);
   ID = find(bfMark==1|bfMark==4);
   Res(ID) = u(ID,1)-0.0;
\end{verbatim}

%\end{exam}

\begin{figure}[h]
\centering
\tikzset{my node/.code=\ifpgfmatrix\else\tikzset{matrix of nodes}\fi}
{\footnotesize
\begin{tikzpicture}[every node/.style={my node},scale=0.2]
\draw[thick] (0,0)--(15,0)--(15,4)--(30,4)--(30,0)--(60,0)--(60,16)--(30,16)--(30,12)--(15,12)--(15,16)--(0,16)--(0,0);

\node[below] at (0,0) {(0,0)\\};
\node[above] at (0,16) {(0,16)\\};
\node[below] at (15,0) {(15,0)\\};
\node[below] at (30,0) {(30,0)\\};
\node[below] at (60,0) {(60,0)\\};
\node[above] at (15,4) {(15,4)\\};
\node[above] at (30,4) {(30,4)\\};

\node[above] at (7.5,-0.5) {$\raisebox{.5pt}{\textcircled{\raisebox{-.9pt} {1}}}$\\};
\node[below] at (7.5,16.5) {$\raisebox{.5pt}{\textcircled{\raisebox{-.9pt} {1}}}$\\};
\node[above] at (45,-0.5) {$\raisebox{.5pt}{\textcircled{\raisebox{-.9pt} {1}}}$\\};
\node[below] at (45,16.5) {$\raisebox{.5pt}{\textcircled{\raisebox{-.9pt} {1}}}$\\};
\node[left] at (60.5,8) {$\raisebox{.5pt}{\textcircled{\raisebox{-.9pt} {1}}}$\\};

\node[left] at (15.5,2) {$\raisebox{.5pt}{\textcircled{\raisebox{-.9pt} {3}}}$\\};
\node[left] at (15.5,14) {$\raisebox{.5pt}{\textcircled{\raisebox{-.9pt} {3}}}$\\};
\node[right] at (29.5,2) {$\raisebox{.5pt}{\textcircled{\raisebox{-.9pt} {3}}}$\\};
\node[right] at (29.5,14) {$\raisebox{.5pt}{\textcircled{\raisebox{-.9pt} {3}}}$\\};
\node[above] at (22.5,3.5) {$\raisebox{.5pt}{\textcircled{\raisebox{-.9pt} {3}}}$\\};
\node[below] at (22.5,12.5) {$\raisebox{.5pt}{\textcircled{\raisebox{-.9pt} {3}}}$\\};

\node[right] at (-0.5,8) {$\raisebox{.5pt}{\textcircled{\raisebox{-.9pt} {2}}}$\\};
\end{tikzpicture}
}
\caption{The domain $\Omega$ and the marking of the boundary segments for Example~\ref{exam-combustion2d}.}
\label{fig:combustion_domain}
\end{figure}
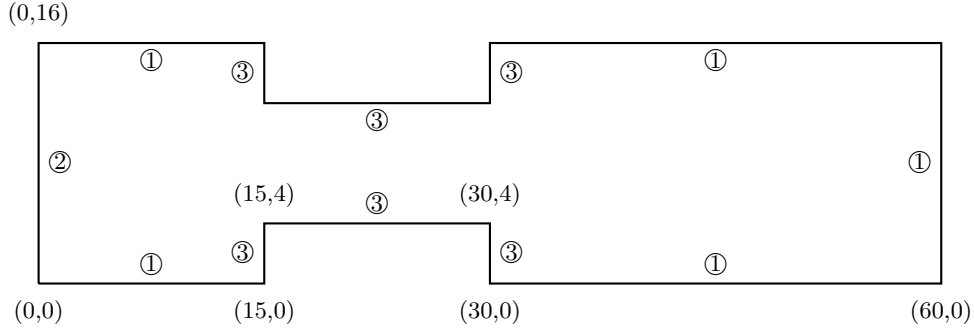

% section 3.3
\subsection{A combustion model in 2D}
% exam: combustion model in 2D
%\begin{exam}
\label{exam-combustion2d}
This example, implemented in \verb|ex2d_combustion.m|, is the IBVP for a combustion model
(a system of two PDEs) in 2D (see \cite{Lang1998}),
\begin{equation}
\label{combustion2d-1}
\begin{cases}
\theta_t = \theta_{xx} + \theta_{yy} + \frac{\beta^2}{2 Le} Y e^{-\frac{\beta (1-\theta)}{(1-\alpha (1-\theta))}},
&\quad (x,y) \in \Omega , \quad t \in (0, 60] \\
Y_t = \frac{1}{Le} (Y_{xx} + Y_{yy}) - \frac{\beta^2}{2 Le} Y e^{-\frac{\beta (1-\theta)}{(1-\alpha (1-\theta))}},
&\quad (x,y) \in \Omega , \quad t \in (0, 60]
\end{cases}
\end{equation}
subject to the boundary conditions
\begin{equation}
\label{combustion2d-1-BC}
\begin{cases}
\theta = 1, \quad Y = 0, & \text{ on bfMark = 2}\\
\frac{\partial \theta}{\partial n} = 0, \quad \frac{\partial Y}{\partial n} = 0, & \text{ on bfMark = 1}\\
\frac{\partial \theta}{\partial n} = -k \theta , \quad \frac{\partial Y}{\partial n} = 0, & \text{ on bfMark = 3}\\
\end{cases}
\end{equation}
and the initial condition
\begin{equation}
\label{combustion2d-1-IC}
\begin{cases}
\theta= 1, \quad Y = 0, & \text{ for } x \le 7.5 \\
\theta = e^{7.5-x}, \quad Y = 1- e^{Le (7.5-x)}, & \text{ for } x > 7.5
\end{cases} 
\end{equation}
where $\Omega$ is shown and the boundary segments are marked as in Fig.~\ref{fig:combustion_domain}
and $Le = 1$, $\alpha = 0.8$, $\beta = 10$, and $k=0.1$.
The analytical expression of the exact solution is not available. The weak formulation reads as
\begin{align}
\label{combustion2d-2}
& \int_\Omega \left ( \theta_t v^{(1)} + \theta_x v^{(1)}_x + \theta_y v^{(1)}_y 
- \frac{\beta^2}{2 Le} Y v^{(1)} e^{-\frac{\beta (1-\theta)}{(1-\alpha (1-\theta))}} \right ) d x dy
+ \int_{\Gamma_{\text{bfMark = 3}}} k \theta v^{(1)} d s \\
& \qquad + \int_\Omega \left (Y_t v^{(2)} + \frac{1}{Le} Y_x v^{(2)}_x + \frac{1}{Le} Y_y v^{(2)}_y 
+ \frac{\beta^2}{2 Le} Y v^{(2)} e^{-\frac{\beta (1-\theta)}{(1-\alpha (1-\theta))}}\right ) d x dy
\notag \\
& \qquad \qquad \qquad \qquad \qquad \qquad. \qquad \qquad \qquad
= 0,\quad \forall v^{(1)}, v^{(2)} \in V
\notag
\end{align}
where $V = \{ v \in H^1(\Omega), v = 0 \text{ on bfMark = 2}\}$.
The definition of this example in the code is given as
\begin{verbatim}
   % define PDE system and BCs

   pdedef.bfMark = ones(Nbf,1);  
   Xbfm = (X(tri_bf(:,1),:)+X(tri_bf(:,2),:))*0.5;
   pdedef.bfMark(Xbfm(:,1) < 1e-8) = 2;
   pdedef.bfMark(abs(Xbfm(:,1)-15) < 1e-8) = 3;
   pdedef.bfMark(abs(Xbfm(:,1)-30) < 1e-8) = 3;
   pdedef.bfMark((abs(Xbfm(:,2)-4) < 1e-8) & ...
       (Xbfm(:,1) > 15 & Xbfm(:,1) < 30)) = 3; 
   pdedef.bfMark((abs(Xbfm(:,2)-12) < 1e-8) & ...
       (Xbfm(:,1) > 15 & Xbfm(:,1) < 30)) = 3; 

   pdedef.bftype = ones(Nbf,npde);
   pdedef.bftype(pdedef.bfMark==1|pdedef.bfMark==3,:) = 0; 
   pdedef.volumeInt = @pdedef_volumeInt;
   pdedef.boundaryInt = @pdedef_boundaryInt;
   pdedef.dirichletRes = @pdedef_dirichletRes;
   
... ...

function F = pdedef_volumeInt(du, u, ut, dv, v, x, t, ipde)
    
    beta = 10;
    alpha = 0.8;
    Le = 1;
    w = beta^2/(2*Le)*u(:,2).*exp(-beta*(1-u(:,1))./(1-alpha*(1-u(:,1))));
    if (ipde==1)
       F = ut(:,1).*v+du(:,1).*dv(:,1)+du(:,2).*dv(:,2) - w.*v;
    else
       F = ut(:,2).*v+(du(:,3).*dv(:,1)+du(:,4).*dv(:,2))/Le + w.*v;
    end 

function G = pdedef_boundaryInt(du, u, v, x, t, ipde, bfMark)

   k = 0.1;
   G = zeros(size(x,1),1);
   if ipde==1
      ID = find(bfMark==3);
      G(ID) = k*u(ID,1).*v(ID);
   end

function Res = pdedef_dirichletRes(u, x, t, ipde, bfMark)

   Res = zeros(size(x,1),1);
   ID = find(bfMark==2);
   if (ipde==1)
      Res(ID) = u(ID,1)-1;
   else
      Res(ID) = u(ID,2)-0;
   end
 \end{verbatim}
%\end{exam}

% section 3.4
\subsection{Poisson's equation in 3D}
% exam: poisson equation in 3D
%\begin{exam}
\label{exam-poisson3d}
This example, implemented in \verb|ex3d_poisson.m|, is the BVP for Poisson's equation in 3D,
\begin{equation}
\label{poisson3d-1}
- (u_{xx} + u_{yy} + u_{zz}) = f,
\quad (x,y,z) \in \Omega \equiv (0,1)\times (0, 1) \times (0,1)
\end{equation}
subject to the boundary conditions
\begin{equation}
\label{poisson3d-1-BC}
\begin{cases}
\frac{\partial u}{\partial x} = 2\pi\sin(3 \pi y)\sin(\pi z), & \quad (x,y,z) \text{ on } \Gamma_N \\
u = u_{exact}(x,y,z), & \quad (x,y,z) \text{ on } \Gamma_D
\end{cases}
\end{equation}
where $\Gamma_N = \{ x = 1, 0< y<1, 0 < z < 1\}$, $\Gamma_D = \partial \Omega \setminus \Gamma_N$, and
$f$ is chosen such that the exact solution of this example is 
\begin{equation}
\label{poisson3d-1-exac}
u_{exact}(x, y, z) = \sin(2\pi x) \sin(3 \pi y) \sin(\pi z) .
\end{equation}
The weak formulation of this example reads as
\begin{align}
\label{poisson3d-2}
\int_\Omega ( u_x v_x + u_y v_y + u_z v_z - f v ) d x dy d z
+ \int_{\Gamma_N} (- 2\pi\sin(3 \pi y)\sin(\pi z)) v(1,y, z) dy d z  = 0,\quad \forall v \in V
\notag
\end{align}
where $V = \{ w \in H^1(\Omega), w = 0 \text{ on } \Gamma_D \}$.
The definition of this example in the code is given as
\begin{verbatim}
   % define PDE system and BCs

   pdedef.bfMark = ones(Nbf,1);
   Xbfm = (X(tri_bf(:,1),:)+X(tri_bf(:,2),:)+X(tri_bf(:,3),:))/3;
   pdedef.bfMark(Xbfm(:,1)>1-1e-8) = 2; % for x=1
   pdedef.bftype = ones(Nbf,npde);
   pdedef.bftype(pdedef.bfMark==2,npde) = 0; % neumann bc for x=1
   
   pdedef.volumeInt = @pdedef_volumeInt;
   pdedef.boundaryInt = @pdedef_boundaryInt;
   pdedef.dirichletRes = @pdedef_dirichletRes;
 
 ... ...
 
function F = pdedef_volumeInt(du, u, ut, dv, v, x, t, ipde)
    F = 14*pi^2*sin(2*pi*x(:,1)).*sin(3*pi*x(:,2)).*sin(pi*x(:,3));
    F = du(:,1).*dv(:,1)+du(:,2).*dv(:,2)+du(:,3).*dv(:,3)-F.*v(:); 

function G = pdedef_boundaryInt(du, u, v, x, t, ipde, bfMark)
   G = zeros(size(x,1),1);
   ID = find(bfMark==2);
   G(ID) = -2*pi*sin(3*pi*x(ID,2)).*sin(pi*x(ID,3)).*v(ID);

function Res = pdedef_dirichletRes(u, x, t, ipde, bfMark)
   Res = u(:,1) - uexact(t,x);
\end{verbatim}
%\end{exam}

%---------------------------------------------------
%\bibliographystyle{plain}
%\bibliographystyle{abbrv}
%\bibliography{/Users/wzhuang/tex/bib/mmesh}

\end{document}